\documentclass[10pt,twoside,english]{myamsart}
\advance\oddsidemargin by -1.0cm
\advance\evensidemargin by -1.0cm
\advance\topmargin by -1.0cm 
\usepackage{amssymb}
\usepackage{babel}
\usepackage{amstext}
\usepackage{amscd}   
\usepackage{epsfig}  
\usepackage{rotating}
\let\mathg\mathfrak
\theoremstyle{plain}

\newtheorem{lem}{Lemma}[section]
\newtheorem{thm}{Theorem}[section]
\newtheorem{prop}{Proposition}[section]
\theoremstyle{definition}
\newtheorem{exa}{Example}[section]
\newtheorem{NB}{Remark}[section]

\renewcommand{\labelenumi}{$(\theenumi)$}
\newcommand{\bdm}{\begin{displaymath}}
\newcommand{\edm}{\end{displaymath}}
\newcommand{\be}{\begin{equation}}
\newcommand{\ee}{\end{equation}}
\newcommand{\ba}[1]{\begin{array}{#1}}
\newcommand{\ea}{\end{array}}

\newcommand{\btab}{\begin{tabular}}
\newcommand{\etab}{\end{tabular}}

\newcommand{\C}{\ensuremath{\mathbb{C}}}
\newcommand{\R}{\ensuremath{\mathbb{R}}}

\newcommand{\Z}{\ensuremath{\mathbb{Z}}}

\newcommand{\vrho}{\ensuremath{\varrho}}
\newcommand{\SU}{\ensuremath{\mathrm{SU}}}
\newcommand{\U}{\ensuremath{\mathrm{U}}}

\newcommand{\SO}{\ensuremath{\mathrm{SO}}}
\newcommand{\Spin}{\ensuremath{\mathrm{Spin}}}

\begin{document}
\def\haken{\mathbin{\hbox to 6pt{%
                 \vrule height0.4pt width5pt depth0pt
                 \kern-.4pt
                 \vrule height6pt width0.4pt depth0pt\hss}}}
    \let \hook\intprod
\setcounter{equation}{0}
%
%------ draw title page -----
%
\thispagestyle{empty}
%
%\hbox to \hsize{%
%  \vtop{} \hfill
%  \vtop{\hbox{PRELIMINARY VERSION}}}
%------------------------------
\date{\today}
%----------------------------------------------------------
\title{On Types of non-integrable geometries}
%----------------------------------------------------------
%
% author and address
%
%-------------------------------------------
%
\author{Thomas Friedrich}
%-------------------------------------------
\address{\hspace{-5mm} 
Thomas Friedrich\newline
Institut f\"ur Mathematik \newline
Humboldt-Universit\"at zu Berlin\newline
Sitz: WBC Adlershof\newline
D-10099 Berlin, Germany\newline
{\normalfont\ttfamily friedric@mathematik.hu-berlin.de}}
%
%-----------------------------------------------------------
\thanks{Supported by the SFB 288 "Differential geometry
and quantum physics" of the DFG}
%-----------------------------------------------------------
\subjclass{Primary 53 C 25; Secondary 81 T 30}
%-----------------------------------------------------------
%-----------------------------------------------------------
\begin{abstract}
%---------------
We study the types of non-integrable $\mathrm{G}$-structures on Riemannian 
manifolds. In particular, geometric 
types admitting a connection with totally skew-symmetric torsion are characterized. $8$-dimensional manifolds equipped with a $\Spin(7)$-structure play
a special role.
Any geometry of that type admits a unique connection with totally
skew-symmetric torsion. Under weak conditions on the structure group
we prove that this geometry is the only one with this property. Finally, we 
discuss the automorphism group of a Riemannian manifold with a fixed
non-integrable $\mathrm{G}$-structure.
\end{abstract}
%-------------
\maketitle
%----------------
\tableofcontents
%----------------
\pagestyle{headings}
%
%
%-------------- body of the document ------------------------------------------
%
%---------------------------------------------------------------------------- 
\section{Introduction}\noindent
%----------------------------------------------------------------------------
Riemannian manifolds equipped with additional geometric structures occur in
many situations and have interesting properties. The most important
structures are almost complex structures and almost contact metric structures. Moreover, in special
dimensions we have exceptional geometries resulting from the list of
exceptional Lie groups, for example there is a $7$-dimensional representation
of the group $\mbox{G}_2$ and a $26$-dimensional representation of the group
$\mbox{F}_4$. In case the Riemannian geometry is compatible
with the additional, geometric structure we call it integrable. The compatibility condition means that the geometric structure under consideration is 
parallel with respect to the Levi-Civita connection or, equivalently, the
Riemannian holonomy group reduces to the subgroup preserving the geometric
structure. Examples are K\"ahler manifolds, Calabi-Yau manifolds, parallel
$\mbox{G}_2$-structures in dimension $7$, parallel $\Spin(7)$-structures in 
dimension $8$ and symmetric spaces. However, there
are many interesting Riemannian manifolds equipped with non-integrable
geometric structures. This happens in any case for almost contact metric
structures in odd dimensions, there are (non-K\"ahler) hermitean manifolds in
even dimensions and non-symmetric, homogeneous spaces. Usually the Riemannian holonomy group of these manifolds is the full orthogonal group.
Consequently, they are of general type in the sense of holonomy 
theory and cannot be distinguished from this point of view.\\

\noindent
A lot of work has been done in order to understand special
non-integrable geometries. In case the geometric structure can be defined
by some tensor $\mathcal{T}$, one considers its Riemannian covariant derivative
$\nabla^{LC}\mathcal{T}$. It is a $1$-form with values in the representation 
space of the tensor. The decomposition of the corresponding tensor product 
under
the action of the group $\mbox{G}$ preserving the tensor $\mathcal{T}$ 
yields the different
classes of non-integrable geometric structures. For any class of
non-integrable geometries one derives a differential equation characterizing
the class and involving the tensor $\mathcal{T}$. This program was developed,
for example, for almost hermitean manifolds (Gray/Hervella 
\cite{GrayHerv}), for $\mbox{G}_2$-structures in dimension $7$ 
(Fernandez/Gray \cite{FerGray}), 
for $\Spin(7)$-structures in dimension $8$ (Fernandez \cite{Fer}) 
and for almost contact metric structures (Chinea/Gonzales \cite{ChiGo}). \\

\noindent
Some years ago I became interested in $16$-dimensional Riemannian
manifolds with a $\Spin(9)$-structure (see \cite{Fri1}, \cite{Fri2}). There
the situation is slightly different, since a structure of that type is
not defined by a single tensor. Therefore, I looked for another method
in order to introduce a classification of non-integrable 
$\mbox{G}$-structures. 
The theory of principal fibre bundles and connections  yields the idea that a 
classification of non-integrable $\mbox{G}$-structures can be based on the difference $\Gamma$ between the Levi-Civita connection and the canonical 
$\mbox{G}$-connection 
induces on the $\mbox{G}$-structure. In some sense $\Gamma$ measures the 
non-integrability of the $\mbox{G}$-structure in a natural way. It is a $1$-form 
defined on the manifold with values in the subspace $\mathg{m}$ orthogonal to 
the Lie algebra 
$\mathg{g}$. At the same time A. Swann 
(see \cite{Swann}, \cite{CleySwa}) and A. Fino (see \cite{Fino}) considered this $1$-form 
for different reasons, too (see also Chiossi/Salamon in \cite{ChiSa}). \\

\noindent
Let us define the different classes of 
non-integrable $\mbox{G}$-structures as the irreducible components of the
representation $\R^n \otimes \mathg{m}$. If the geometric
structure is given by a tensor, this point of view is completely equivalent
to the approach described before. One of the aims
of this note is to explain that one obtains all the known 
results in a unified way. The approach seems to be a kind of ``folklore'' for 
some people,
but even in differential geometry it is not as popular as it should be. 
It will turn out that the reproduction of some classical results 
cited before becomes much less computational in our approach.
It also has the advantage of being applicable to geometric
structures not defined by a tensor. For example, we discuss 
irreducible $\SO(3)$-structures on $5$-dimensional manifolds, $\Spin(9)$-structures on $16$-dimensional manifolds as well as $\mbox{F}_4$-structures on
$26$-dimensional manifolds. Some problems concerning non-integrable geometric
structures can be solved  immediately from this point of view. In string 
theory 
one wants to know which types of geometric structures 
admit affine connections $\nabla$ with totally skew-symmetric torsion 
(see \cite{Gauntlett}, \cite{Strominger}). 
It turns out that the answer depends mainly on the decomposition of two 
representations into irreducible components. An interesting example are 
$8$-dimensional Riemannian manifolds with a $\Spin(7)$-structure. It is well 
known (see \cite{Iv}) that \emph{any} $\Spin(7)$-structure admits a unique 
connection with totally skew-symmetric
torsion. In this paper we prove that under certain weak conditions on the
structure group this is the only geometry with this property. Finally, we
study the automorphism group
of non-integrable geometric structures.
%---------------------------------------------------------------------------- 
\section{$\mbox{G}$-structures on Riemannian manifolds}\noindent
%----------------------------------------------------------------------------
%
Let $\mbox{G} \subset \SO(n)$ be a closed subgroup of the orthogonal group
and decompose the Lie algebra
\bdm
\mathg{so}(n) \ = \ \mathg{g} \oplus \mathg{m} 
\edm
into the subalgebra $\mathg{g}$ and its orthogonal complement $\mathg{m}$. 
We denote by $\mbox{pr}_{\mathg{g}}$ and $\mbox{pr}_{\mathg{m}}$ the 
projections of the Lie algebra $\mathg{so}(n)$ onto $\mathg{g}$ and 
$\mathg{m}$, respectively. 
Consider an oriented Riemannian manifold $(M^n, g)$ and denote 
its frame bundle by 
$\mathcal{F}(M^n)$. It is a principal $\SO(n)$-bundle 
over $M^n$.
A $\mbox{G}$\emph{-structure} of $M^n$ is a reduction 
$\mathcal{R} \subset \mathcal{F}(M^{n})$ of the frame bundle to the subgroup 
$\mbox{G}$. The Levi-Civita connection is a $1$-form $Z$ on 
$\mathcal{F}(M^{n})$ with values in the Lie algebra $\mathg{so}(n)$.
We restrict the Levi-Civita connection to $\mathcal{R}$ and decompose 
it with respect to the decomposition of the Lie algebra $\mathg{so}(n)$:
\bdm
Z\big|_{T(\mathcal{R})} \ := \ Z^* \, \oplus \ \Gamma \, .
\edm
Then, $Z^*$ is a connection in the principal $\mbox{G}$-bundle $\mathcal{R}$
and $\Gamma$ is a tensorial $1$-form of type Ad, i.\,e., a $1$-form on $M^{n}$ 
with values in the associated bundle $\mathcal{R} \times_{\mbox{G}} 
\mathg{m}$. The triple $(M^n, g, \mathcal{R})$ is an \emph{integrable}
$\mbox{G}$-structure if the $1$-form $\Gamma$ vanishes, i.\,e., the Levi-Civita
connection preserves the $\mbox{G}$-structure $\mathcal{R}$. Many interesting
geometric structures are not of that type. In this paper we consider mainly 
\emph{non-integrable} geometric structures, $\Gamma \neq 0$.  We
introduce a general classification of these structures using the
$\mbox{G}$-type of the $1$-form $\Gamma$. More precisely, the 
$\mbox{G}$-representation $\R^{n} \otimes \mathg{m}$ splits into irreducible 
components. The different \emph{non-integrable types} of $\mbox{G}$-structures 
are defined -- via the decomposition of $\Gamma$ -- as the irreducible 
$\mbox{G}$-components of the representation $\R^{n} \otimes \mathg{m}$. 
Let us give a local formula for $\Gamma$. Fix an orthonormal frame
$e_1, \ldots , e_n$ adapted to the reduction $\mathcal{R}$. The connection
forms $\omega_{ij} := g(\nabla^{LC} e_i , e_j)$ of the Levi-Civita connection
define a $1$-form $\Omega : = (\omega_{ij})$ with values in the Lie algebra
$\mathg{so}(n)$ of all skew-symmetric matrices. The form $\Gamma$ is the
$\mathg{m}$-projection of $\Omega$,
\bdm
\Gamma \ = \ \mbox{pr}_{\mathg{m}}(\Omega) \ = \ 
\mbox{pr}_{\mathg{m}}(\omega_{ij}) \ .
\edm\\
\noindent
The case that the subgroup $\mbox{G}$ is the isotropy group
of some tensor $\mathcal{T}$ is of special interest. Suppose that there is a 
faithful representation $\vrho : \SO(n) \rightarrow 
\SO(V)$ and a tensor $\mathcal{T} \in V$ such that
\bdm
\mbox{G} \ = \ \big\{ g \in \SO(n) : \vrho(g)\mathcal{T} = \mathcal{T}\big\} \ .
\edm
Then a 
$\mbox{G}$-structure is a triple $(M^n, g, \mathcal{T})$ consisting
of a Riemannian manifold equipped with an additional tensor field. The
Riemannian covariant derivative is given by the formula
\bdm
\nabla^{LC} \mathcal{T} \ = \ \vrho_*(\Gamma)(\mathcal{T}) \, ,
\edm
where $\vrho_* : \mathg{so}(n) \rightarrow \mathg{so}(V)$ is the differential
of the representation. $\nabla^{LC}\mathcal{T}$ is an element
of $\R^n \otimes V$. The algebraic $\mbox{G}$-types of $\nabla^{LC}\mathcal{T}$
define the algebraic $\mbox{G}$-types of $\Gamma$ and vice versa. Indeed, 
we have
\begin{prop}
The $\mathrm{G}$-map
\bdm
\R^n \otimes \mathg{m} \longrightarrow \R^n \otimes \mathrm{End}(V) \longrightarrow \R^n \otimes V
\edm 
given by $\Gamma \rightarrow \rho_*(\Gamma)(\mathcal{T})$ is injective.
\end{prop}
\begin{proof}  
If $\rho_*(\Gamma)(\mathcal{T}) = 0$, then the endomorphism 
$\rho_*(\Gamma(X))$ stabilizes $\mathcal{T}$ for any vector $X \in \R^n$, 
i.e., $\rho_*(\Gamma(X)) \in \rho_*(\mathg{g})$. Since the representation is 
faithful, we conclude that $\Gamma(X) \in \mathg{g}$. On the other hand, 
we have 
$\Gamma(X) \in \mathg{m}$, i.e., $\Gamma \equiv 0$.
\end{proof}
\noindent
The covariant derivative $\nabla^{LC}\mathcal{T}$ 
has been
used for the classification of geometric structures -- see the examples. 
The approach presented here uses the 
$1$-form $\Gamma$ and applies even in case that the geometric
structure is not defined by a tensor. Moreover, in many situations it
is simpler to handle the $\mbox{G}$-type of $\Gamma$ then the $\mbox{G}$-type 
of the covariant derivative.\\
\begin{prop} 
%---------
If the group $\mathrm{G}$ does not coincide with the full group
$\SO(n)$, then the $\mathrm{G}$-representation $\R^n$ is always one of the
components of the representation $\R^n \otimes \mathg{m}$. 
\end{prop}
\begin{proof} 
Indeed, consider the map
\bdm
\R^n \longrightarrow \R^n \otimes \mathg{m} , \quad \quad X \longrightarrow 
\sum_{i=1}^n e_i \otimes \mbox{pr}_{\mathg{m}}(e_i \wedge X) \ ,
\edm
where $e_1, \ldots , e_n$ is an orthonormal basis in $\R^n$. Suppose
that a vector $X$ belongs to its kernel. Then, for any vector $Y$, the 
exterior product $X \wedge Y$ is an element of the Lie algebra $\mathg{g}$. 
Since the commutator of two elements again belongs to the Lie algebra 
$\mathg{g}$, we conclude that the exterior product $Y \wedge Z$ of two vectors 
orthogonal to the vector $X$ is in $\mathg{g}$, i.e., 
$\mathg{g} = \mathg{so}(n)$. 
\end{proof}
\noindent
Geometrically this fact reflects the conformal 
transformation of a $\mbox{G}$-structure. Let $(M^n,g,\mathcal{R})$ be a 
Riemannian manifold with a fixed geometric
structure and denote by $\hat{g} := e^{2f} \cdot g$ a conformal transformation
of the metric. There is a natural identification of the frame bundles
\bdm
\mathcal{F}(M^n, g) \ \cong \mathcal{F}(\hat{M}^n, \hat{g})
\edm
and a corresponding $\mbox{G}$-structure $\mathcal{\hat{R}}$. On the 
infinitesimal level, the conformal change is defined by the $1$-form 
$df$.
\subsection{$\SO(3)$-structures in dimension $5$}
%-----------------------------------------------------
%
The group $\SO(3)$ has a unique, real, irreducible
representation in dimension $5$. We consider the corresponding non-standard 
embedding
$\SO(3) \subset \SO(5)$ as well as the decomposition
\bdm
\mathg{so}(5) \ = \ \mathg{so}(3) \oplus \mathg{m}^7 \ .
\edm
It is well known that the $\SO(3)$-representation $\mathg{m}^7$ is the
unique, real, irreducible representation in dimension $7$. We decompose the
tensor product into irreducible components
\bdm
\R^5 \otimes \mathg{m}^7 \ = \ \R^3 \oplus \R^5 \oplus \mathg{m}^7 \oplus
\mbox{E}^{9} \oplus \mbox{E}^{11} \ .
\edm
There are five basic types of $\SO(3)$-structures on $5$-dimensional Riemannian
manifolds. The symmetric space $\SU(3)/\SO(3)$ is an example of a $5$-dimensional Riemannian manifold with an integrable $\SO(3)$-structure, ($\Gamma = 0$).
\subsection{Almost complex structures in dimension $6$}
%-----------------------------------------------------
%
Let us consider $6$-dimensional Riemannian manifolds $(M^6, g, 
\mathcal{J})$ with an almost complex structure $\mathcal{J}$. The subgroup
$\U(3) \subset \SO(6)$ describes a geometric structure of that type. We
decompose the Lie algebra
\bdm
\mathg{so}(6) \ = \ \mathg{u}(3) \oplus \mathg{m}
\edm
and remark that the $\U(3)$-representation in $\R^6$ is the real representation underlying $\Lambda^{1,0}$ and, similarly, $\mathg{m}$ is the real 
representation underlying $\Lambda^{2,0}$. We decompose the complexification
under the action of $\U(3)$:
\bdm
\Big(\R^6 \otimes \mathg{m}\Big)^{\C} \ = \ \Big(\Lambda^{1,0} \otimes \Lambda^{2,0}
\oplus \Lambda^{1,0} \otimes \Lambda^{0,2}\Big)_{\R}^{\C} \ .
\edm
The symbol $( \ldots )_{\R}^{\C}$ means that we understand the complex
representation as a real representation and complexify it. Next we split
the complex $\U(3)$-representations
\bdm
\Lambda^{1,0} \otimes \Lambda^{2,0} \ = \ \C^3 \otimes \Lambda^2(\C^3) \ 
= \ \Lambda^{3,0} \oplus \mbox{E}^8 \ ,
\edm
\bdm
\Lambda^{1,0} \otimes \Lambda^{0,2} \ = \ \C^3 \otimes 
\Lambda^2(\overline{\C}^3) \ = \ \C^3 \otimes \Lambda^2(\C^3)^* \ = \
(\C^3)^* \oplus \mbox{E}^6 \ .
\edm
$\mbox{E}^6$ and $\mbox{E}^8$ are irreducible $\U(3)$-representations
of complex dimension $6$ and $8$, respectively. Finally we obtain
\bdm
\R^6 \otimes \mathg{m} \ = \ \Lambda^{3,0} \oplus (\C^3)^* \oplus \mbox{E}^6 
\oplus \mbox{E}^8 \ .
\edm
Consequently, $\R^6 \otimes \mathg{m}$ splits into four irreducible representations of real dimensions $2, \, 6, \, 12$ and $16$, i.e., there are four basic types of $\U(3)$-structures on
$6$-dimensional Riemannian manifolds (Gray/Hervella-classification - 
see \cite{GrayHerv}). In case we restrict the structure group to $\SU(3)$, we 
obtain two trivial summands in the decomposition of $\R^6 \otimes \mathg{su}(3)^{\perp}$ corresponding to \emph{nearly K\"ahler} manifolds (see 
\cite{Gray},\cite{Gray2},\cite{Swann}). Almost hermitean manifolds of that 
type 
have  special properties in real dimension $n=6$. They are Einstein manifolds 
(see \cite{Gray2}), the differential equation describing the nearly
K\"ahler manifolds is
\bdm
\big(\nabla^{LC}_X\mathcal{J}\big)\big(X\big) \ = \ 0, \quad \quad 
\nabla^{LC}\mathcal{J} \ \neq \ 0
\edm
and, finally, these are precisely the $6$-dimensional manifolds with real
Killing spinors (see \cite{Grun}).
\subsection{$\mathrm{G}_2$-structures in dimension $7$}
%------------------------------------------------------
%
We consider $7$-dimensional Riemannian manifolds equip\-ped with a 
$\mbox{G}_2$-structure. Since the group $\mbox{G}_2$ is the isotropy
group of a $3$-form $\omega^3$ of general type, a $\mbox{G}_2$-structure
is a triple $(M^7, g, \omega^3)$ consisting of a $7$-dimensional Riemannian
manifold and a $3$-form $\omega^3$ of general type at any point.  We decompose the $\mbox{G}_2$-representation (see \cite{Friedrich&I1})
\bdm
\R^7 \otimes \mathg{m} \ = \ \R^1 \oplus \R^7 \oplus \Lambda^2_{14} \oplus 
\Lambda^3_{27}  
\edm
and, consequently, there are four basic types of non-integrable 
$\mbox{G}_2$-structure.  In this way we obtain the Fernandez/Gray-classification of $\mbox{G}_2$-structures (see \cite{FerGray}). The different
types of $\mbox{G}_2$-structures can be characterized by differential
equations. For example, a $\mbox{G}_2$-structure is of type $\R^1$ 
(\emph{nearly parallel} structures) if and only if there exists a number $\lambda$ such
that
\bdm
 d \omega^3 \ = \ \lambda \cdot (* \omega^3)
\edm
holds. Again, in dimension $n=7$ this condition is equivalent to the existence of a real
Killing spinor  (see \cite{FKMS}). The $\mbox{G}_2$-structures of type $\R^1 \oplus 
\Lambda^3_{27}$ (\emph{cocalibrated} structures) are characterized by the
condition that the $3$-form is coclosed, $\delta \omega^3 = 0$. In general, 
the differential equations
for any type of $\mbox{G}_2$-structure involving the $3$-form $\omega^3$ 
were derived in \cite{FerGray}. In the spirit of the approach
of this paper one can find the computations in \cite{Friedrich&I1}.
\subsection{$\Spin(7)$-structures in dimension $8$}
%--------------------------------------------------
Let us consider $\Spin(7)$-structures on $8$-dimensional Riemannian
manifolds. The subgroup $\Spin(7) \subset \SO(8)$ is the real 
$\Spin(7)$-representation $\Delta_7 = \R^8$. The complement 
$\mathg{m} = \R^7$ is the standard $7$-dimensional representation and the 
$\Spin(7)$-structures on an $8$-dimensional Riemannian manifold $M^8$ 
correspond to the irreducible components of the tensor product
\bdm
\R^8 \otimes \mathg{m} \ = \ \R^8 \otimes \R^7 \ = \ \Delta_7 \otimes \R^7 \ = \ \Delta_7 \oplus \mbox{K}  \ = \ \R^8 \oplus \mbox{K} \, ,
\edm
where $\mbox{K}$ denotes the kernel of the Clifford multiplication 
$\Delta_7 \otimes \R^7 \to \Delta_7$. It is well known that 
$\mbox{K}$ is an irreducible $\Spin$-representation, i.e., there are two 
basic types of $\Spin(7)$-structures (the Fernandez-classification of
$\Spin(7)$-structures - see \cite{Fer}).
\subsection{$\Spin(9)$-structures in dimension $16$}
%---------------------------------------------------
The group $\Spin(9)$ is an interesting subgroup of $\SO(16)$. The
representation in $\R^{16}$ is irreducible. 
We consider $16$-dimensional Riemannian manifolds with $\Spin(9)$-structures. 
Again, we split the $\Spin(9)$-representation into four irreducible
components
\bdm
\R^{16} \otimes \mathg{m} \ = \ \R^{16} \oplus \mathcal{P}_1(\R^9) \oplus 
\mathcal{P}_2(\R^9) \oplus \mathcal{P}_3(\R^9)\ . 
\edm
Consequently, there are four basic types of non-integrable 
$\Spin(9)$-structures
on $16$-dimensional Riemannian manifolds (see \cite{Fri1}, \cite{Fri2}). 
A $\Spin(9)$-structure on a $16$-dimensional Riemannian manifold has an
associated $8$-form (but is \emph{not} defined by this $8$-form in the sense 
we have explained above). Some
differential equations for the different types of $\Spin(9)$-structures
involving the $8$-form have been computed in the paper \cite{Fri1}.
\subsection{$\mbox{F}_4$-structures in dimension $26$}
%---------------------------------------------------
We consider the subgroup $\mbox{F}_4 \subset \SO(26)$ and $26$-dimen\-sional
Riemannian manifolds with a $\mbox{F}_4$-structure. The orthogonal complement
\bdm
\mathg{so}(26) \ =  \ \mathg{f}_4 \oplus \mathg{m}^{273}
\edm
is the unique, irreducible, $273$-dimensional representation of $\mbox{F}_4$.
We compute the decomposition
\bdm
\R^{26} \otimes \mathg{m}^{273} \ = \ \R^{26} \oplus \mathg{f}_4 \oplus
\mathg{m}^{273} \oplus \mbox{E}^{324} \oplus \mbox{E}^{1053} 
\oplus \mbox{E}^{1274} \oplus \mbox{E}^{4096} \ .
\edm
Consequently, there are seven basic types of $\mbox{F}_4$-structures
in dimension $26$.

\section{\mbox{G}-connections with totally skew-symmetric torsion} 
%-----------------------------------------------------------------
%
\noindent
An interesting question in some models in string theory (see \cite{Friedrich&I1}) is to ask of which geometric structures admit a connection $\nabla$ preserving the structure and with totally skew-symmetric torsion. In order to 
formulate the general condition, let us introduce the maps
\bdm
\Theta_1 : \Lambda^3({\Bbb R}^n) \longrightarrow {\Bbb R}^{n} \otimes \mathg{m}, \quad
\Theta_2 : \Lambda^3({\Bbb R}^n) \longrightarrow {\Bbb R}^{n} \otimes \mathg{g}
\edm 
given by the formulas
\bdm
\Theta_1(T) \ := \ \sum_i (\sigma_i \haken T) \otimes \sigma_i, 
\quad
\Theta_2(T) \ := \ \sum_j (\mu_j \haken T) \otimes \mu_j \ ,
\edm
where $\sigma_i$ is an orthonormal basis in $\mathg{m}$ and $\mu_j$ is an
orthonormal basis in $\mathg{g}$. Then we have
\begin{thm}\label{thm2} (see \cite{Fri2})
%----------------------
A $\mathrm{G}$-structure $\mathcal{R} \subset \mathcal{F}(M^{n})$ of a 
Riemannian manifold admits a connection $\nabla$ with totally skew-symmetric 
torsion if and only if the $1$-form $\Gamma$ belongs to the image 
of $\Theta_1$, 
\bdm
2 \cdot \Gamma \ = \ - \, \Theta_1(T) \ . 
\edm
In this case the $3$-form $T$ is the torsion form of the connection.
\end{thm}
\begin{proof} Suppose there exists a connection $\nabla$ with totally skew-symmetric torsion $T$. We compare it with the Levi-Civita connection and obtain
the relation
\bdm
\nabla_X Y \ = \ \nabla^{LC}_X Y + \frac{1}{2} \cdot T(X,Y, \ast) \ .
\edm
Moreover, the definition of the $1$-form $\Gamma$ as well as the $\mbox{G}$-connection $Z^*$ yield the equation
\bdm
\nabla^{LC}_X Y \ = \ \nabla^{Z^*}_X Y + \Gamma(X)Y \ .
\edm
Finally, since $\nabla$ preserves the $\mbox{G}$-structure, there
exists a $1$-form $\beta$ with values in the Lie algebra $\mathg{g}$ such that
\bdm
\nabla_X Y \ = \ \nabla^{Z^*}_X Y + \beta(X)Y \ .
\edm
Combining the three formulas we obtain, for any vector $X$, the equation
\bdm
2 \cdot \beta(X) \ = \ 2 \cdot \Gamma(X) + T(X, \ast, \ast) \ .
\edm
We project onto the subspace $\mathg{m}$. Since $\beta(X)$ belongs
to the Lie algebra $\mathg{g}$, we conclude that $\Gamma$ should be in the 
image
of $\Theta_1$,  $\Theta_1(T) = - 2 \cdot \Gamma$.
\end{proof} 

\noindent
Theorem \ref{thm2} only decides
which geometric types admit connections with totally skew-symmetric torsion.
However, if the geometric structure is defined by a tensor $\mathcal{T}$, 
one prefers to express the torsion form $T$ of the connection $\nabla$ 
directly by this tensor $\mathcal{T}$. Formulas of that type were 
computed for almost complex structures,
almost contact metric structures, $\mbox{G}_2$-structures and $\Spin(7)$-structures (see \cite{Friedrich&I1}, \cite{Friedrich&I2}, \cite{Friedrich&I3}, 
\cite{Iv} and \cite{Agricola}). 
\begin{exa} We consider $5$-dimensional Riemannian manifolds with an
$\SO(3)$-structure. Then we obtain
\bdm
\R^5 \otimes \mathg{m}^7 \ = \ \R^3 \oplus \R^5 \oplus \mathg{m}^7 \oplus
\mbox{E}^{9} \oplus \mbox{E}^{11}, \quad \Lambda^3(\R^5) \ = \ \R^3 \oplus 
\mathg{m}^7 \ .
\edm
In particular, a conformal change of an $\SO(3)$-structure does \emph{not} 
preserve the property that the structure admits a connection with totally 
skew-symmetric torsion.
\end{exa}
\begin{exa} In the case of almost complex structures in dimension $6$, we have
\bdm
\R^6 \otimes \mathg{m} \ = \ \Lambda^{3,0} \oplus (\C^3)^* \oplus \mbox{E}^6 
\oplus \mbox{E}^8, \quad \Lambda^3(\R^6) \ = \ \Lambda^{3,0} \oplus (\C^3)^* \oplus \mbox{E}^6 \ .
\edm
Consequently, an almost complex manifold $(M^6, g, \mathcal{J})$ admits
a connection with totally skew-symmetric torsion if and only if the
$\mbox{E}^8$-part of $\Gamma$ vanishes (see \cite{Friedrich&I1}). In case
the connection exists, it is unique.
\end{exa}
\begin{exa} 
%----------
In dimension $7$ we decompose the $\mbox{G}_2$-representation
(see \cite{Friedrich&I1})
\bdm
\Lambda^3(\R^7) \ = \ \R^1 \oplus \R^7 \oplus \Lambda^3_{27}, \quad 
\R^7 \otimes \mathg{m} \ = \ \R^1 \oplus \R^7 \oplus \Lambda^2_{14} \oplus 
\Lambda^3_{27} \ . 
\edm
Consequently, a $\mbox{G}_2$-structure admits a connection with totally
skew-symmetric torsion if and only if it is of type $\R^1 \oplus \R^7 
\oplus \Lambda^3_{27}$. These condition describes the conformal changes
of cocalibrated $\mbox{G}_2$-structures. In case the connection exists, 
it is unique.
\end{exa}
\begin{exa}\label{exa8}  
%----------
Let us consider $\Spin(7)$-structures on $8$-dimensional Riemannian
manifolds. Here we find
\bdm
\R^8 \otimes \mathg{m} \ = \ \Delta_7 \oplus \mbox{K}, \quad \quad 
\Lambda^3(\R^8) \ = \ \Delta_7 \oplus \mbox{K} \ ,
\edm
i.\,e.,  $\Lambda^3(\R^8) \to \R^8 \otimes \mathg{m}$ is an isomorphism.
Theorem \ref{thm2} yields immediately that \emph{any $\Spin(7)$-structure on an
$8$-dimensional Riemannian manifold admits a unique connection with totally skew-symmetric torsion} (see \cite{Iv}).
\end{exa}
\begin{exa} 
%----------
In case of $\mbox{G} = \Spin(9)$, we have
\bdm
\R^{16} \otimes \mathg{m} \ = \ \R^{16} \oplus \Lambda^3(\R^{16}) \oplus 
\mathcal{P}_3(\R^9) \ ,
\edm
and the $\R^{16}$-component is \emph{not} contained in $\Lambda^3(\R^{16}) = 
 \mathcal{P}_1(\R^9) \oplus \mathcal{P}_2(\R^9)$. 
A conformal change of a $\Spin(9)$-structure does \emph{not} 
preserve the property that the structure admits a connection with totally 
skew-symmetric torsion. 
\end{exa}
\begin{exa} In dimension $26$ and for the subgroup $\mbox{F}_4 \subset \SO(26)$
we have
\bdm
\R^{26} \otimes \mathg{m}^{273} \ = \ \R^{26} \oplus \mathg{f}_4 \oplus
\mathg{m}^{273} \oplus \mbox{E}^{324} \oplus \mbox{E}^{1053} 
\oplus \mbox{E}^{1274} \oplus \mbox{E}^{4096}, \quad 
\Lambda^3(\R^{26}) \ = \ \mathg{m}^{273} \oplus \mbox{E}^{1053} \oplus 
\mbox{E}^{1274} \ .
\edm
In particular, a conformal change of an $\mbox{F}_4$-structure does \emph{not} 
preserve the property that the structure admits a connection with totally 
skew-symmetric torsion. 
\end{exa}
\section{The automorphism group of non-integrable \mbox{G}-structures} 
\noindent
We consider a Riemannian manifold $(M^n, g, \mathcal{R})$ with a fixed
geometric structure. Since the Lie algebra $\Lambda^2(\R^n) = \mathg{so}(n)
= \mathg{g} \oplus \mathg{m}$ splits, the bundle of $2$-form $\Lambda^2(M^n)$
decomposes into two subbundles. They are associated with the reduction 
$\mathcal{R}$ of the frame bundle and we denote these two bundles again by $\mathg{g}$ and $\mathg{m}$, respectively,
\bdm
\Lambda^2(M^n) \ = \ \mathg{g} \oplus \mathg{m} \ .
\edm
Let $X$ be a Killing vector field. Then the covariant derivative 
$\nabla^{LC}X \in \Gamma(T \otimes T)$ is skew-symmetric. In fact, if
we understand $X$ as a $1$-form on the manifold, the covariant derivative
of $X$ coincides with the exterior differential, 
\bdm
\nabla^{LC} X \ = \ \frac{1}{2} \cdot dX \ .
\edm
We suppose now that the $\mbox{G}$-structure $\mathcal{R}$ admits a unique
connection $\nabla$ with totally skew-symmetric torsion $T$. Then 
$\nabla X \in \Gamma(T \otimes T)$ is a skew-symmetric tensor, too. Moreover,
we have
\bdm
\nabla X \ = \ \nabla^{LC} X - \frac{1}{2} \cdot (X \haken T) \ .
\edm
\begin{thm}\label{thm3} Let $(M^n, g, \mathcal{R})$ be a $\mathrm{G}$-structure
and suppose that there exists a unique $\mathrm{G}$-connection $\nabla$ with
totally skew-symmetric torsion $T$. If a Killing vector field $X$ is an 
infinitesimal transformation of the $\mathrm{G}$-structure, then
\bdm
\mathcal{L}_X T \ = \ 0\ , \quad \ [X,  \nabla_Y Z] - \nabla_Y [X,Z] \ = \ 
\nabla_{[X,Y]} Z \ .  
\edm
The $2$-form $\nabla X \in \mathg{g}$ belongs to the subbundle $\mathg{g}$.
In particular, we have 
\bdm
\mbox{pr}_{\mathg{m}}(dX) \ = \ \mbox{pr}_{\mathg{m}}(X \haken T) \ .
\edm
\end{thm}
\begin{proof} Since $\nabla$ is the unique $\mbox{G}$-connection with totally
skew-symmetric torsion, any transformation of $\mathcal{R}$ should
preserve the connection and its torsion form, i.e., the first two conditions
are necessary. In fact, the condition $\nabla X \in
\mathg{g}$ characterizes the infinitesimal transformations
preserving a $\mbox{G}$-structure. Let us - for completeness - give 
the argument. The covariant 
derivative of a vector field with respect to an affine 
metric connection can be computed via the formula
\bdm
g (\nabla_Y X, \, Z)(p) \ := \ \frac{d}{dt} \, g(df_t(p)(Y), 
\, \tau^{\nabla}_t(Z) ) \ ,
\edm
where $f_t : M^n \rightarrow M^n$ is the $1$-parameter group generated by
the vector field $X$ and $\tau^{\nabla}_t$ denotes the parallel displacement
along the curve $f_t(p)$. Fix a basis $e_1, \ldots , e_n \in \mathcal{R}_p$
in the $\mbox{G}$-structure at the point $p \in M^n$ and denote by $A_{ij}(t)$ the matrix defined by
\bdm
df_t(e_i) \ := \ \sum_{j=1}^n A_{ij}(t) \cdot \tau^{\nabla}_t(e_j) \ .
\edm
The endomorphism $\nabla X(p)$ is given by the matrix $\big(A_{ij}'(0)\big)$. 
If the $1$-parameter group $f_t$ preserves the structure $\mathcal{R}$, then
the matrix $\big(A_{ij}(t)\big)$ belongs to the subgroup $\mbox{G}$, i.e., 
$\nabla X$ is a $2$-form in $\mathg{g}$.
\end{proof}
\begin{NB}
The formula $\mbox{pr}_{\mathg{m}}(dX) \ = \ \mbox{pr}_{\mathg{m}}(X 
\haken T)$ was derived in case of a nearly parallel 
$\mbox{G}_2$-structure in \cite[Theorem 6.2]{FKMS} (notice that there
is a sign error). Indeed, a
nearly parallel $\mbox{G}_2$-structure admits a unique connection 
with totally skew-symmetric torsion $T$, which was computed in
\cite{Friedrich&I1}, Example 5.2. Using these expression for $T$ we obtain 
from Theorem \ref{thm3} the formula of Theorem 6.2 in \cite{FKMS}.
\end{NB}
\noindent
The invariance of the torsion form restricts the dimension of the 
automorphism $\mathcal{G}(\mathcal{R})$. Denote by $\mbox{G}_T$ the isotropy group
of $T \in \Lambda^3(\R^n)$ and $dT \in \Lambda^4(\R^n)$ inside of $\mbox{G}$. 
Then we have
\bdm
\mbox{dim}\big(\mathcal{G}(\mathcal{R})\big) \ \leq \ n + \mbox{dim} (\mbox{G}_T) \ .
\edm
The group $\mbox{G}_T$ preserves the Ricci tensor of the unique connection
$\nabla$ as well as the symmetric endomorphism $T_{imn} \cdot T_{jmn}$. 
These geometric objects have been computed in several cases and can be used 
in the computation of the isotropy group of the torsion form.
\begin{exa} Denote by $\mbox{H}^6$ the $6$-dimensional Heisenberg group. 
There exists a left-invariant, cocalibrated $\mbox{G}_2$-structure $\omega^3$ 
on the $7$-dimensional Lie group $\mbox{H}^6 \times \R^1$. In 
\cite{Friedrich&I1} we computed its torsion form $T$:
\bdm
\omega^3 \ = \ e_{127} + e_{135} - e_{146} - e_{236} - e_{245} + 
e_{347} + e_{567}\, , \quad 
T \ = \ e_5 \wedge (e_{13} - e_{67}) + e_4 \wedge (e_{37} + e_{16}) \ .
\edm
Moreover, the Ricci tensor $\mbox{Ric}^{\nabla}$ of the unique connection
with totally skew-symmetric torsion as well as the symmetric endomorphism
$T_{imn} \cdot T_{jmn}$ are given by the formulas
\bdm
\mbox{Ric}^{\nabla} \ = \ \mbox{diag} (-2 , \,  0 ,\, -2 ,\, 0 ,\, 0 ,\, -2 
,\, -2) \, , \quad T_{imn} \cdot T_{jmn} \ = \ \mbox{diag} ( 4 , \,  0 ,\, 4 ,\, 4 ,\, 4 ,\, 4 ,\, 4) \ .
\edm
A transformation preserving the geometric structure preserves the Ricci 
tensor $\mbox{Ric}^{\nabla}$ and the symmetric form $T_{imn} \cdot T_{jmn}$, 
too. Consequently, for the Lie algebra
$\mathg{g}_T$ of the group $\mbox{G}_T$ we obtain the necessary conditions
$ \omega_{2\alpha} =  0, \ \omega_{4\beta} = 0, \ \omega_{5\gamma} = 0$ 
for any $1 \leq \alpha \leq 7$, $\beta \neq 5$ and $\gamma \neq 4$. Combining these $14$ equations with the equations defining the Lie algebra $\mathg{g}_2$ 
inside of $\mathg{so}(7)$ (see \cite{FKMS} or \cite{Friedrich&I1})
we obtain seven nontrivial parameters $\omega_{13}, \omega_{16}, 
\omega_{17}, \omega_{36}, \omega_{37}, 
\omega_{67}, \omega_{45}$, related by three equations
\bdm
\omega_{13} \ = \ - \, \omega_{67} \, , \quad \omega_{16} \ = \ \omega_{37} \, , 
\quad \omega_{17} + \omega_{36} + \omega_{45} \ = \ 0 \ .
\edm 
We understand the skew-symmetric matrix $\Omega := (\omega_{ij})$ as a vector
field on $\R^7$ and compute the Lie derivative $\mathcal{L}_{\Omega}T$ of the
torsion form,
\bdm
\mathcal{L}_{\Omega}T \ = \ 2 \cdot \omega_{13} \cdot (e_{147} - e_{346}) 
+ 2 \cdot \omega_{16} \cdot (e_{356} - e_{157}) + 2 \cdot \omega_{17} 
\cdot ( e_{357} + e_{156} + e_{134} - e_{467}) \ .
\edm
Consequently, the Lie group $\mbox{G}_T$ is one-dimensional and its Lie algebra
is described by two parameters $\omega_{36}, \, \omega_{45}$ and one equation:
\bdm
\omega_{36} + \omega_{45} \ = \ 0 \ .
\edm
\end{exa}
\begin{exa} The product $\mbox{N}^6 \times \R^1$ of the $3$-dimensional,
complex, solvable Lie group $\mbox{N}^6$ by $\R^1$ admits a left invariant
cocalibrated $\mbox{G}_2$-structure. In \cite{Friedrich&I1} we derived its
torsion,
\bdm
T \ = \ 2 \cdot (e_{256} - e_{234}) \ .
\edm
A straightforward calculation yields that the subgroup $\mbox{G}_T \subset 
\mbox{G}_2$ is a maximal torus of $\mbox{G}_2$. A basis of its Lie algebra is 
given by the following two matrices
\bdm
\left[\begin{array}{rrrrrrr} 0&0&0&0&0&0&-2\\
0&0&0&0&0&0&0\\ 0&0&0&0&0&1&0\\ 0&0&0&0&1&0&0\\ 0&0&0&-1&0&0&0\\
0&0&-1&0&0&0&0\\ 2&0&0&0&0&0&0\end{array}\right],\quad
\left[\begin{array}{rrrrrrr} 0&0&0&0&0&0&0\\
0&0&0&0&0&0&0\\ 0&0&0&1&0&0&0\\ 0&0&-1&0&0&0&0\\ 0&0&0&0&0&-1&0\\
0&0&0&0&1&0&0\\ 0&0&0&0&0&0&0\end{array}\right]\,.
\edm
\end{exa}
\section{A characterization of \mbox{Spin}(7)-structures} 
%---------------------------------------------------------
\noindent
Let us once again return to Example \ref{exa8}. An $8$-dimensional Riemannian
manifold equipped with a $\Spin(7)$-structure admits a unique connection 
preserving the $\Spin(7)$-structure with totally skew-symmetric torsion 
(see \cite{Iv}). In general, fix  a compact, connected subgroup 
$\mbox{G} \subset \SO(n)$ and consider $\mbox{G}$-geometries. A
Riemannian $\mbox{G}$-manifold $(M^n, g, \mathcal{R})$ admits 
a unique
$\mbox{G}$-connection with totally skew-symmetric torsion if and only if the 
$\mbox{G}$-representations
\bdm
\Theta_1 : \Lambda^3({\Bbb R}^n) \longrightarrow {\Bbb R}^{n} \otimes \mathg{m}
\edm
are isomorphic (see Theorem \ref{thm2}). We will prove that under a
certain condition on the group $\mbox{G}$ only
the case of $n = 8$ and $\mbox{G} = \Spin(7)$ is possible.
\begin{lem}\label{lem}  
Let $\mathrm{G}$ be a compact, connected Lie group and denote by
$\mathrm{T}$ its maximal torus. Then the inequality
\bdm
\mathrm{dim(G)} \ \le \ 4 \cdot \big(\mathrm{dim(T)}\big)^2
\edm
holds. Moreover, if no exceptional Lie algebra occurs in the decomposition of 
the Lie algebra $\mathg{g}$, then we have
\bdm
\mathrm{dim(G)} \ \le \ 3 \cdot \big(\mathrm{dim(T)}\big)^2 .
\edm
\end{lem} 
\begin{proof} 
Remark that the inequality holds for the product of 
two groups $\mbox{T}_1 \subset \mbox{G}_1 , \ \mbox{T}_2 \subset \mbox{G}_2$ 
in case it holds for $\mbox{G}_1$ and  $\mbox{G}_2$. Indeed, $\mbox{T}_1 
\times \mbox{T}_2$ is a maximal torus in $\mbox{G}_1 \times \mbox{G}_2$ and we obtain
\bdm
\mathrm{dim(G_1 \times G_2)} \ \le \ 4 \cdot \big( (\mathrm{dim(T_1)}^2 + 
(\mathrm{dim(T_2)}^2\big) \ \le \ 4 \cdot \big(\mathrm{dim(T_1)} + 
\mathrm{dim(T_2)} \big)^2 \ .
\edm
We split the Lie algebra $\mathg{g} = \mathg{g}_1 \oplus \cdots \oplus 
\mathg{g}_l \oplus \mathg{z}$ into
the simple ideals $\mathg{g}_i$ and its center $\mathg{z}$. 
Unless $\mathg{g}_i$ is a classical simple Lie algebra, we know that its 
dimension is bounded by $3 \cdot (\mathrm{dim(T}_i))^2$ (see \cite{FH}). 
For the exceptional Lie algebras we obtain
\bdm
\mathrm{dim(G_2)} \ = \ 14, \quad 4 \cdot (\mathrm{dim(T)})^2 \ = \ 16; \quad 
\quad \mathrm{dim(F_4)} \ = \ 52, \quad 4 \cdot (\mathrm{dim(T)})^2 \ = 
\ 64  ;
\edm
\bdm
\ \ \  \mathrm{dim(E_6)} \ = \ 72, \quad 4 \cdot (\mathrm{dim(T)})^2 \ 
= \ 144;
\quad  \mathrm{dim(E_7)} \ = \ 133, \quad 4 \cdot (\mathrm{dim(T)})^2 \ 
= \ 196; 
\edm
and $\mathrm{dim(E_8)} = 248, \ 4 \cdot (\mathrm{dim(T)})^2 = 256$.
\end{proof}
\noindent
Let $h \in \mbox{G}$ be an element of the Lie group and denote by
\bdm
\mathrm{Z}(h) \ := \ \big\{g \in \mathrm{G} : g \cdot h = h \cdot g \big\}, \quad \mathrm{z} \ := \ \mathrm{dim}\, \mathrm{Z}(h) 
\edm 
its centralizer as well as its dimension. We agree to say that a subgroup 
$G \subset \SO(n)$ of dimension $\mbox{g} := \mbox{dim}\mbox(G)$ has the 
\emph{involution property} if one of the following conditions is satisfied:
\begin{enumerate}
\item
$n^2 \ \not= \ 3 \cdot \mathrm{g} + 1 \ .$ 
\item
$n^2 \ = \ 3 \cdot \mathrm{g} + 1$, but there does not exist a pair 
$(h,p)$ consisting
of an involution $h \in \mbox{G}$ and an even integer 
$0 < p < n$ such that
\bdm
3 \cdot (\mathrm{g} - \mathrm{z}) \ = \ 
2 \cdot p \cdot \big(\sqrt{3 \cdot \mathrm{g} + 1} - p) \ .
\edm
\end{enumerate}
Solving the latter equation with respect to $p$ we obtain
\bdm
p \ = \ \frac{1}{2}\big(\sqrt{1 + 3\cdot \mathrm{g}} \, \pm \,  
\sqrt{6 \cdot \mathrm{z} - 3 \cdot \mathrm{g} + 1}\big) \, .
\edm
\begin{NB} Using the representation of $\mbox{G}$ in $\R^n$ we can formulate 
this condition in a more 
geometric way. Fix an involution 
$h \in \mbox{G} \subset \SO(n)$ and consider the two symmetric spaces
\bdm
\frac{\mathrm{G}}{\mathrm{Z}(h)} \ \subset \ 
\frac{\SO(n)}{\mathrm{Z}_{\SO(n)}(h)} \ = \ \mathrm{G}_{n,p} \ ,
\edm
where $\mbox{G}_{n,p}$ denotes the Grassmannian manifold of all oriented
$p$-planes in $\R^n$ (p even). By the involution property we want to 
exclude the case that $n \ = \ \sqrt{3 \cdot \mathrm{g} + 1}$  
and the ratio
of the dimensions of these two symmetric spaces is $\frac{2}{3}$. 
\end{NB}
\begin{exa}
%--------------
Consider the group $\SU(3)$. Then $\sqrt{3 \cdot \mathrm{g} + 1} = 5$ and 
the dimension of the centralizer
of the involution $h = \mbox{diag}(1,-1,-1)$ equals $\mbox{z} = 4$. 
In particular,
$6 \cdot \mbox{z} - 3 \cdot \mbox{g} +1 = 1$ and $p = 2, 3$. 
Nevertheless there does not exist
a subgroup of $\SO(5)$ that is isomorphic to $\SU(3)$, i.e., $\SU(3)$ has the
involution property. 
\end{exa}
\begin{exa} 
%----------
In case of $\Spin(7)$, we have $\sqrt{3 \cdot \mathrm{g} + 1} = 8$, 
but there is no involution $h \in \Spin(7)$ such that
\bdm
3 \cdot \big( 21 - \mathrm{dim}\ \mathrm{Z}(h) \big) \ = \  2 \cdot p \cdot 
\big( 8 - p \big)
\edm
for an even number $p$. More generally, we have
\end{exa}
\begin{prop} Any compact simple Lie group 
$\mathrm{G}$ has the involution property. 
\end{prop}
\begin{proof} 
For the exceptional Lie groups $\mbox{G}_2, \mbox{F}_4, \mbox{E}_6, 
\mbox{E}_7, \mbox{E}_8$ the number
$\sqrt{3 \cdot \mbox{g} + 1}$ is not an integer. 
For the classical groups the irreducible symmetric spaces 
$\mbox{G}/\mbox{Z}(h)$ given by an involution $h \in \mbox{G}$ are well known
(see \cite{Goodman}, chapter 11.2.4)
\bdm
\SU(m)/\mathrm{S}\big(\U(r) \times \U(m-r)\big), \quad 
\SO(m)/\SO(r) \times \SO(m-r), \quad \mathrm{Sp}(m)/\mathrm{Sp}(r) 
\times \mathrm{Sp}(m-r) \ .
\edm
The dimension of the group $\mbox{G} = \SU(m)$ equals $(m^2 -1)$ and 
we obtain the restriction $n = \sqrt{3 \cdot m^2 - 2} \approx \sqrt{3} \cdot m$. On the
other hand, the lowest \emph{real} dimension of an $\SU(m)$-representation 
is $2 \cdot m, \, (m \geq 3)$. This means that even
in case $n = \sqrt{3 \cdot m^2 - 2}$ is an integer, there is no subgroup
of $\SO(n)$ that is isomorphic to $\SU(m)$. A similar argument applies to the group
$\mbox{Sp}(m)$. Finally, we discuss the case that $\mbox{G}$ is locally
isomorphic to $\SO(m)$. Then $n = \sqrt{\frac{3}{2} \cdot m \cdot (m-1) +1}
\approx \sqrt{\frac{3}{2}} \cdot m$. Taking into account the dimensions
of all irreducible real $\SO(m)$-representations we conclude that the
embedding $\SO(m) \rightarrow \SO(n)$ is the standard inclusion
\bdm
h \ \longrightarrow \ \left(\!\!\begin{array}{cc}h&0\\ 0& I_{n-m}\end{array}\!\!
\right)\, ,
\edm
and the embedding of the symmetric spaces is the usual inclusion of two
Grassmannian manifolds
\bdm
\frac{\mathrm{G}}{\mathrm{Z}(h)} \ = \  \mathrm{G}_{m,p}\  \longrightarrow \
\frac{\SO(n)}{\mathrm{Z}_{\SO(n)}(h)} \ = \ \mathrm{G}_{n,p} \ .
\edm
In particular, $p$ is bounded by $m$, $p < m$. The dimension condition 
yields $3 \cdot m - p = 2 \cdot n$ and
together with $2 \cdot n^2 = 3 \cdot m \cdot (m-1) +2$ we see that there is no solution $(p, m, n)$ of these two equations satisfying the condition $p < m$. 
\end{proof}
\begin{thm} 
%----------
Let $\mathrm{G} \subset \SO(n)$ be a compact, connected group with
the involution property. Decompose the Lie algebra into 
$\mathg{so}(n) = \mathg{g} \oplus \mathg{m}$ and suppose that the 
$\mathrm{G}$-representations
\bdm
\Theta_1 : \Lambda^3({\Bbb R}^n) \longrightarrow {\Bbb R}^{n} \otimes \mathg{m}
\edm
are isomorphic. Then $n = 8$, $\mathrm{G} = \Spin(7)$ and the representation
is the unique irreducible representation of $\Spin(7)$ in $\R^8$.
\end{thm}
\begin{proof}
%------------
Denote by $\chi , \chi^* : \mbox{G} \rightarrow \R^1$ the characters of the
$\mbox{G}$-representation in $\R^n$ and of the adjoint representation
$\mathg{g}$. Since $\Lambda^3(\R^n)$ is isomorphic to $\R^n \otimes \mathg{m}$ by assumption, 
an easy computation yields the following functional equation between these two
characters
\bdm
3 \cdot \chi(h) \cdot \chi^*(h) \ = \ \chi^3(h) - \chi(h^3),  \quad \ h \in 
\mathrm{G}
\edm
(see \cite{FH}, page $381$). Evaluating the characters at the element $h = e$ we obtain a formula relating the dimensions of the group $\mbox{G}$ to the dimension of the representation $\R^n$, 
\bdm
n^2 \ = \ 3 \cdot \mathrm{dim(G)} + 1 \ .
\edm 
Fix a maximal torus $\mbox{T}^t \subset \mbox{G}$ of dimension $t$ and 
denote
by $h_1:= e, \, h_2, \cdots , \, h_{2^t}$ its elements of order two,
\bdm
h_i^2 \ = \ e, \quad h_i \cdot h_j \ = \ h_j \cdot h_i \ .
\edm
The character equation simplifies for each of these elements, 
\bdm
3 \cdot \chi(h_i) \cdot \chi^*(h_i) \ = \ \chi^3(h_i) - \chi(h_i) \ = \ 
\chi(h_i) \cdot \big(\chi^2(h_i) - 1 \big) \ .
\edm
Suppose that $\chi(h_{i_0}) \not= 0$ for some index $2 \leq i_0 \leq 2^t$. Then we obtain
$3 \cdot \chi^*(h_{i_0}) \ = \ \chi^2(h_{i_0}) - 1$.
Using the equations
\bdm
\chi^*(h_{i_0}) \ = \ 2 \cdot \mathrm{dim\, Z}(h_{i_0}) - 
\mathrm{dim(G)}, 
\quad n^2 \ = \ 3 \cdot \mathrm{dim(G)} + 1, \quad \chi(h_{i_0}) \ 
= \ n - 4 \cdot q
\edm
we see that the latter equation contradicts the involution property of the group $\mbox{G}$ except for $q = n/2$ and $h_{i_0} = - \, \mbox{Id}_{\R^n}$. Consequently, the character $\chi(h_i) = 0$ 
vanishes for any element $h_i \not= \pm \, \mbox{Id}_{\R^n}$. Then 
the number
\bdm
k \ := \ \frac{n}{2^{t-1}} \ \in \ \Z 
\edm
must be an integer. Indeed, denote by $\mbox{H}$ the finite group consisting
of all involutions $h_i$. Consider the space  
$(\R^n)^{\mathrm{H}}$ of all $\mbox{H}$-invariant vectors and its dimension
(see \cite{FH}, page 16),
\bdm
\mathrm{dim}\, (\R^n)^{\mathrm{H}} \ = \ \frac{1}{2^t} \sum_{i=1}^{2^t} 
\chi(h_i) \ .
\edm 
If the involution $(- \, \mbox{Id}_{\R^n})
\not\in \mbox{H}$ does not belong to the subgroup, we obtain
\bdm
\mathrm{dim}\, (\R^n)^{\mathrm{H}} \ = \ \frac{n}{2^t} \ .
\edm 
If $(- \, \mbox{Id}_{\R^n})$ is an element of $\mbox{H}$, 
we choose
a subgroup $\mbox{H}_0 \subset \mbox{H}$ of order two not containing this
involution. In this case we obtain
\bdm
\mathrm{dim}\, (\R^n)^{\mathrm{H_0}} \ = \ \frac{n}{2^{t-1}} \ .
\edm
By the previous Lemma \ref{lem} we have
the inequality
\bdm
4 \cdot t^2 \ \geq \ \mathrm{dim(G)} \ = \ \frac{1}{3}(4^{t-1}\cdot k^2 -1)
\ \geq \ \frac{1}{3}(4^{t-1} - 1) \ .
\edm
In particular, the rank of the compact group $\mbox{G}$ is bounded by five, $t \leq 5$. The cases $ t = 1, 2$ or $4$ can be directly excluded by
the conditions
\bdm
3 \cdot \mathrm{dim(G)} + 1\ = \ 4^{t-1} \cdot k^2, \quad \mathrm{dim(G)} \ 
\leq 4 \cdot t^2 \ .
\edm 
Let us discuss the case of $t = 3$. Then we obtain the conditions
\bdm
\mathrm{dim(G)} \ \leq \ 36, \quad 3 \cdot \mathrm{dim(G)} + 1 \ = \ 16 
\cdot k^2 \ .
\edm
If $k = 1$, the dimension of the group equals five, $\mbox{dim(G)} = 5$, and
the dimension $n$ of the real representation is given by the formula 
$n^2 = 3 \cdot \mbox{dim(G)} +1 = 16$. Therefore, 
the group $\mbox{G}$ is a compact subgroup of rank three in
$\SO(4)$, a contradiction. In case of $k=2$ we obtain $\mbox{dim(G)} = 21$ and
$n=8$. The group is a $21$-dimensional subgroup of rank $3$ in $\SO(8)$,
i.e., $\mbox{T}^3 \subset \mbox{G} = \Spin(7)\subset \SO(8)$. The cases
that $k \geq 3$ are impossible $(t = 3)$.\\

\noindent
Finally, we discuss the case of $t=5$. Then we obtain the conditions
\bdm
\mathrm{dim(G)} \ \leq \ 100, \quad 3 \cdot \mathrm{dim(G)} + 1 \ = \ 256 
\cdot k^2 \ .
\edm
The parameters $k \geq 2$ are impossible and $k=1$ yields an $85$-dimensional
subgroup $\mbox{G} \subset \SO(16)$ of rank five. Since $85 = \mbox{dim(G)} 
\not\leq 3 \cdot t^2$, by Lemma \ref{lem} the decomposition of the Lie
algebra $\mathg{g}$ into simple Lie algebras must contain one of the
exceptional algebras $\mathg{g}_2$ or $\mathg{f}_4$. Again, this is impossible. Suppose, for example, that $\mathg{g} = \mathg{g}_2 \oplus 
\mathg{g}^*$. Then the compact Lie algebra $\mathg{g}^*$ has dimension
$71$ and rank $3$ and these parameters contradict Lemma \ref{lem}. A similar
argument excludes the second exceptional summand.
\end{proof}
%
% 
%------------------------------------------
%\addcontentsline{toc}{section}{Literature}
%------------------------------------------
    
\end{document}